\documentclass[12pt]{article}

\usepackage{graphicx,enumitem}
\usepackage{graphicx,psfrag,epsf}
{%
\centering\addtocounter{figure}{1}
\begin{enumerate}[%
itemsep=2pt,parsep=0em,
label={(\alph*)},
ref={\thefigure.(\alph*)}
]}%
{\end{enumerate}\addtocounter{figure}{-1}}

\usepackage[margin=1in]{geometry}

\usepackage[utf8]{inputenc}
\usepackage{amsmath}
\usepackage[english]{babel}
\usepackage{amssymb}
\usepackage{amsfonts}
\usepackage{bm}
\usepackage{graphicx}
\usepackage{amsmath}
\usepackage{amsthm}
\usepackage{xcolor}
\usepackage[hidelinks=true]{hyperref}
\usepackage{parskip}
\usepackage{multicol}
\usepackage{float}

\usepackage{verbatim}

\floatstyle{plaintop}
\restylefloat{table}
\usepackage[tableposition=top]{caption}

\usepackage{rotating}
\usepackage{subcaption}
\usepackage{multirow}

\usepackage{chngcntr}
\usepackage{apptools}
\AtAppendix{\counterwithin{lem}{subsection}}
\AtAppendix{\counterwithin{proposition}{subsection}}
\AtAppendix{\counterwithin{equation}{subsection}}

\newtheorem{proposition}{Proposition}[section]

\newtheorem{ejemplo}{Example}[section]

\usepackage{booktabs}
\usepackage{animate}

\usepackage{natbib}
\hypersetup{urlcolor=blue, citecolor=blue, colorlinks=true, linkcolor=blue} 

\numberwithin{equation}{section}
\numberwithin{figure}{section}
\numberwithin{table}{section}

\usepackage{blindtext}
\usepackage{lineno}

\begin{document}

\title{\textbf{Modelling and simulation of multifractal star-shaped particles}}
\author{Alfredo Alegr\'ia \\ \\ Departamento de Matem\'atica, \\ Universidad T\'ecnica Federico Santa Mar\'ia,\\ Valpara\'iso, Chile.}
\maketitle

\begin{abstract}
\noindent The problem of constructing flexible stochastic models to describe the variability in shape of solid particles is challenging. Natural objects often exhibit mono- or multi-fractal features, i.e. irregular shapes and self-similar patterns. This paper presents a general framework for modelling three-dimensional star-shaped particles with a locally variable Hausdorff (or fractal) dimension. In our approach, the radial function of the particle is represented by an anisotropic Gaussian random field on the sphere. We additionally derive a simulation algorithm being parenthetical to the spectral turning bands method proposed in Euclidean spaces. We illustrate the use of our proposal through numerical examples, including a multifractal simulated version of the Earth topography.    
\end{abstract}


\section{Introduction}

Realistic modelling of three-dimensional star-shaped particles has recently attracted great interest in a wide variety of scientific disciplines such as astronomy, material science and medicine, among others \citep{stoyan1994fractals}.  The Hausdorff (or fractal) dimension, a concept coming from the vocabulary of irregular shapes and self-similar patterns, is of paramount importance in the analysis of solid particles because it serves as a powerful mathematical tool to quantify the degree of smoothness or roughness of the surface of the particle. Indeed, recent studies suggest that fractal geometry can be useful for describing the topography of celestial bodies in the solar system \citep{kucinskas1992fractal},  the mechanisms of tumour growth and angiogenesis \citep{sedivy1997fractals} and the morphological features of sand particles \citep{zhou2017three}.

In the last decades, appealing stochastic models for star-shaped objects with constant Hausdorff dimension (also called monofractal objects)  have been proposed by various authors, including \cite{kent2000using}, \cite{hobolth2003spherical}, \cite{ziegel2013stereological} and \cite{hansen2015gaussian}.  Nevertheless, it is not uncommon to observe natural phenomena with multifractal patterns.  For instance,  \cite{gagnon2006multifractal} argue that the Earth topography is multifractal and propose a simulation method in planar domains.  \cite{dellino2002fractal} study the relevance of multifractal  volcanic ash particles from the eruptions of Monte Pilato-Rocche Rosse. The strong evidence of these complex geometries certainly calls for the development of new methodologies.

 In this paper we propose an approach for modelling three-dimensional multifractal star-shaped particles. We follow the flexible model presented by \cite{hansen2015gaussian}, where the radial function of the particle is represented by mean of a Gaussian random field on the sphere.   Specifically, the covariance function of the random field determines the Hausdorff dimension of the surface of the particle. While \cite{hansen2015gaussian} focused on the  study of isotropic random fields, and the mentioned monofractal objects, we extend their approach to the anisotropic case by admitting locally adaptive spatial dependencies. As a result, our model allows for a Hausdorff dimension varying from place to place on the boundary of the particle.  
 
 We additionally derive an algorithm for high resolution simulations, which relies on a spectral representation of the covariance function. We illustrate the use of our proposal through numerical examples, including a simulated version of the Earth topography,  where the discrepancy between the regularities of the continents and the seafloor, as discussed by  \cite{gagnon2006multifractal}, has been taken into account.
 
The paper is organized as follows. In Section \ref{rfs} we review preliminary results about isotropic random fields on the sphere and introduce the definition of the  fractal index. We also describe the modelling of star-shaped random particles in terms of its radial function.  Section \ref{proposal} proposes a general framework for escaping from isotropy, allowing for random fields with locally adaptive regularity properties. The connections between our proposal and a kernel-based method are also discussed. In addition, the associated simulation algorithm is derived.  Section \ref{numerical} illustrates numerical examples, while Section \ref{discusion} concludes the paper with a discussion.


\section{Background}

\label{rfs}

\subsection{Isotropic Gaussian random fields on the sphere}

Let $\mathbb{S}^2  = \{\boldsymbol{x} \in \mathbb{R}^{3}: \boldsymbol{x}^\top \boldsymbol{x}=1\}$ be the two-dimensional unit sphere, where $^\top$ denotes the transpose operator,  and consider a real-valued random field, $\{ Z(\boldsymbol{x}): \boldsymbol{x}\in\mathbb{S}^2\}$, with finite second order moments.  We assume  $Z(\boldsymbol{x})$ to be Gaussian, i.e. for all $k\in\mathbb{N}$ and $\boldsymbol{x}_1,\hdots,\boldsymbol{x}_k\in\mathbb{S}^2$, the random vector $\{Z(\boldsymbol{x}_1),\hdots,Z(\boldsymbol{x}_k)\}^\top$ follows a multivariate Gaussian distribution. Thus, $Z(\boldsymbol{x})$  is completely characterized by its mean function, $\mu(\boldsymbol{x}) = E\{Z(\boldsymbol{x})\}$, $\boldsymbol{x}\in\mathbb{S}^2$, and by its covariance function $C(\boldsymbol{x}_1,\boldsymbol{x}_2) = {\rm cov}\{Z(\boldsymbol{x}_1),Z(\boldsymbol{x}_2)\}$,   $\boldsymbol{x}_1,\boldsymbol{x}_2\in\mathbb{S}^2$.

Let us introduce the geodesic distance on $\mathbb{S}^2$, which is the main ingredient to define the property of \emph{isotropy} of a random field.     For two locations, $\boldsymbol{x}_1$ and $\boldsymbol{x}_2$ in $\mathbb{S}^2$,  their geodesic distance is defined as $d(\boldsymbol{x}_1,\boldsymbol{x}_2) = \arccos\{ \boldsymbol{x}_1^\top \boldsymbol{x}_2 \} \in [0,\pi]$.  We shall equivalently use $d(\boldsymbol{x}_1,\boldsymbol{x}_2)$ or the shortcut $d$ to denote the geodesic distance, whenever no confusion arises. Following \cite{marinucci2011random}, the random field is called (weakly) isotropic if it has a constant mean, and if its covariance function can be written as
\begin{equation}
\label{isotropic-part}
C(\boldsymbol{x}_1, \boldsymbol{x}_2) = K\{ d(\boldsymbol{x}_1,\boldsymbol{x}_2) \}, \qquad \boldsymbol{x}_1, \boldsymbol{x}_2 \in \mathbb{S}^2,
\end{equation} 
for some continuous function $K: [0,\pi]\rightarrow \mathbb{R}$. Thus, the covariance function just depends on the geodesic distance or, equivalently, the inner product.   It is common to call $K$ the \emph{isotropic part} of the covariance function $C$ (see, e.g., \citealp{guella2018unitarily}). For Gaussian random fields, isotropy also implies that the probability distribution of $\{Z(\boldsymbol{x}_1),\hdots,Z(\boldsymbol{x}_k)\}^\top$ is invariant under the group of rotations on $\mathbb{S}^2$ (see \citealp{marinucci2011random}).    

Note that for all $k\in\mathbb{N}$, and for all systems of points $\boldsymbol{x}_1,\hdots, \boldsymbol{x}_k \in \mathbb{S}^2$ and constants $a_1,\hdots,a_k \in \mathbb{R}$, we have that
\begin{equation}
\label{pos_def}
 {\rm var}\left\{  \sum_{i=1}^k a_i Z(\boldsymbol{x}_i)\right\} = \sum_{i=1}^k \sum_{j=1}^k a_i a_j C(\boldsymbol{x}_i,\boldsymbol{x}_j) \geq 0.
\end{equation}
 Condition (\ref{pos_def}), called \emph{semi positive definiteness},  is a necessary and sufficient condition for a covariance function.  In his pioneering paper,  \cite{schoenberg1942} showed that $C$ as in (\ref{isotropic-part}) is semi positive definite if, and only if, its isotropic part $K$ has a series representation in the form
\begin{equation}
\label{schoenberg1}
K(d) = \sum_{n=0}^\infty   b_n   {P}_n(\cos d), \qquad   0 \leq d\leq \pi,
\end{equation}
where ${P}_n$ denotes the Legendre polynomial of degree $n$, and $\{b_n: n\in\mathbb{N}_0\}$ is a sequence of nonnegative coefficients, such that $\sum_{n=0}^\infty b_n < \infty$. Classical inversion formulas yield
$$  b_n = \frac{2n+1}{2} \int_0^\pi P_n(\cos \xi) \sin(\xi) K(\xi) \text{d}\xi, \qquad n\in\mathbb{N}_0.$$
Following terminology introduced by \cite{ziegel2014convolution}, we refer to this sequence as a \emph{Schoenberg sequence}.

The covariance function is often specified to belong to a parametric family whose members are known to be semi positive definite.  For a thorough review of positive definite functions on spheres and a vast list of parametric families, we refer the reader to \cite{gneiting2013strictly}.

\subsection{Fractal index}

The asymptotic behaviour of the isotropic part $K$ near zero, quantified by the fractal index, regulates the degree of smoothness or roughness of the sample paths of the associated isotropic random field. Recent findings in this direction can be found in \cite{hitczenko2012some}, \cite{hansen2015gaussian}, \cite{lang2015isotropic}, \cite{guinness2016isotropic} and \cite{de2018regularity}.  

Formally, a random field with covariance function as in (\ref{isotropic-part}) has fractal index $\alpha>0$ if there exists a constant $c_0 >0$ such that  
\begin{equation}
\label{fractal_index}
K(0)-K(d)  \sim c_0  d^\alpha,
\end{equation}
as  $d\downarrow 0$. The fractal index exists for most parametric families of covariance functions, in which case it is always true that $0 < \alpha \leq 2$, where $\alpha = 2$ and $\alpha \rightarrow 0$ correspond respectively to extreme smoothness and roughness of the sample paths. 

Abelian and Tauberian theorems \citep{bingham1978tauberian,malyarenko2004abelian} relate the asymptotic behaviour of $K$ near zero to that of its Schoenberg sequence near infinity. In consequence, the fractal index can be characterized in terms of the decay of the Schoenberg sequence.   Indeed, a function $f: (0, \infty) \rightarrow (0, \infty)$ is called \emph{slowly varying} at infinity if, for all $r>0$,  
$$ \lim_{t\rightarrow \infty} \frac{f(rt)}{f(t)} = 1.$$ 
 Then, \cite{malyarenko2004abelian} shows that, for $0 < \alpha < 2$,
\begin{equation}
\label{mayarenko2}
K(0)-K(d)  \sim 2^{-\alpha} d^\alpha f(1/d),
\end{equation}
as $d\downarrow 0$, if, and only if,
\begin{equation}
\label{mayarenko}
\sum_{k=n}^\infty b_k \sim f(n) n^{-\alpha},
\end{equation}
as $n\rightarrow \infty$. The implication from (\ref{mayarenko}) to  (\ref{mayarenko2}) is called Abelian, whereas the converse is called Tauberian. The following example plays a fundamental role throughout the manuscript.

\begin{ejemplo}
The Legendre-Mat\'ern covariance function proposed by \cite{guinness2016isotropic} is characterized by the Schoenberg sequence
\begin{equation}
\label{guinness}
b_n = (\tau^2 + n^2)^{-\nu - 1/2}, \qquad n \in\mathbb{N}_0,
\end{equation}
where $\tau$ and $\nu$ are positive parameters. While $\tau$ regulates the practical range (the distance at which the covariance function reaches certain threshold) of the random field, $\nu$ is capable of controlling the smoothness of the sample paths. Indeed, since 
$$ \sum_{k=n}^\infty (\tau^2 + k^2)^{-\nu - 1/2}   \sim n^{-2\nu} (2\nu)^{-1},$$
as $n\rightarrow \infty$, we conclude from (\ref{mayarenko}) that for $0 < \nu < 1$ the associated random field has fractal index $\alpha=2\nu$. This result can be extended to the limit case $\nu=1$, using additional technical arguments in \cite{bingham1978tauberian}.  It is also worth noting that for $\nu>1$ the random field is mean square differentiable, in which case the fractal index is $\alpha=2$, and we refer the reader to \cite{guinness2016isotropic} for details.   

\end{ejemplo}

\subsection{Star-shaped particles}

We represent a particle as a compact set $Y\subset \mathbb{R}^3$, being star-shaped with respect to an interior point $\boldsymbol{o}$, that is, for all $\boldsymbol{x}\in Y$, the line segment from $\boldsymbol{o}$ to $\boldsymbol{x}$ is contained in $Y$.  The set $Y$ is completely characterized by its radial function, defined as  $Z(\boldsymbol{x}) = \max \{ r\geq 0:  \boldsymbol{o}+r\boldsymbol{x} \in Y  \}$,  for $\boldsymbol{x}\in\mathbb{S}^2$.  Accordingly, $Y$ can be represented as
$$  Y = \bigcup_{\boldsymbol{x}\in\mathbb{S}^2} \{\boldsymbol{o} + r\boldsymbol{x} : 0\leq r \leq Z(\boldsymbol{x})\}.$$
Adopting the framework proposed by \cite{hansen2015gaussian}, we model the radial function $Z(\boldsymbol{x})$ as a Gaussian random field on $\mathbb{S}^2$. Since $Z(\boldsymbol{x})$ has potentially negative values, it might be necessary replace $Z(\boldsymbol{x})$ with $Z_c(\boldsymbol{x}) = \max\{c,Z(\boldsymbol{x})\}$ for some $c>0$.  The interior point $\boldsymbol{o}$ can be assumed to be the origin.

 The regularity of the surface of  $Y$ can be mathematically quantified by the Hausdorff dimension \citep{adler1981geometry}. We now turn to a formal definition in terms of ball coverings \citep{hansen2015gaussian}.  For $\epsilon>0$, an $\epsilon$-cover of $Y$ is a countable collection $\{B_i: i=1,2,\hdots\}$ of balls $B_i  \subset \mathbb{R}^3$ of diameter $|B_i|$ less than or equal to $\epsilon$ that covers $Y$. Let
$$  H^\eta(Y) = \lim_{\epsilon\rightarrow 0} \inf \left\{   \sum |B_i|^\eta :   \{B_i: i=1,2,\hdots\} \text{ is an } \epsilon\text{-cover of } Y \right\}$$
be the $\eta$-dimensional Hausdorff measure of $Y$. The Hausdorff dimension of $Y$ is the unique nonnegative number $\eta_0$ such that $H^\eta(Y)  = \infty$ if $\eta < \eta_0$ and  $H^\eta(Y)  = 0$ if $\eta > \eta_0$.

 \cite{hansen2015gaussian} show that in the special case when $Z(\boldsymbol{x})$ is isotropic, the Hausdorff dimension of the particle is $\eta_0 = 3 - \alpha/2$ almost surely, where $\alpha$ is the fractal index of $Z(\boldsymbol{x})$ defined in (\ref{fractal_index}). Hence, it is always true that $2 \leq \eta_0 < 3$. For sets describing traditional smooth shapes, the Hausdorff dimension matches the conventional topological dimension $\eta_0 = 2$.  In contrast, as the Hausdorff dimension exceeds the topological dimension, the set becomes progressively irregular.

We see that when the radial function is characterized by an isotropic random field, the Hausdorff dimension of the surface of the object is constant. Particles with spatially varying Hausdorff dimension shall be obtained by relaxing the hypothesis of isotropy.

\section{Anisotropic Random Fields on the Sphere}
\label{proposal}

\subsection{Locally varying Schoenberg sequences}

 As discussed in the previous section, the construction of versatile models for star-shaped particles relies on the appropriate specification of anisotropic random field models.  We propose an approach to escape from isotropy, based on spatially adaptive Schoenberg sequences.  More precisely, we consider the covariance function 
\begin{equation}
\label{cov_nonstationary}
C(\boldsymbol{x}_1,\boldsymbol{x}_2) = \sum_{n=0}^\infty   \left\{ b_n(\boldsymbol{x}_1)  b_n(\boldsymbol{x}_2) \right\}^{1/2}   {P}_n(\boldsymbol{x}_1^\top \boldsymbol{x}_2),    \qquad \boldsymbol{x}_1, \boldsymbol{x}_2\in\mathbb{S}^2,
\end{equation}
where $\{b_n(\boldsymbol{x}): n\in\mathbb{N}_0\}$ is a sequence of nonnegative functions on $\mathbb{S}^2$, with $\sum_{n=0}^\infty b_n(\boldsymbol{x}) < \infty$, for each $\boldsymbol{x}\in\mathbb{S}^2$.  We call this sequence of functions an \emph{adaptive Schoenberg sequence}.

It is straightforward  to show that (\ref{cov_nonstationary}) yields a semi positive definite function (see  Appendix A for a proper justification). Moreover, (\ref{cov_nonstationary}) is still positive definite if we replace $ \left\{ b_n(\boldsymbol{x}_1)  b_n(\boldsymbol{x}_2) \right\}^{1/2}$ with any positive definite function on $\mathbb{S}^2\times \mathbb{S}^2$. However, we shall see that (\ref{cov_nonstationary}) is sufficiently general to achieve models with flexible fractal properties. 

This construction is convenient because when both $\boldsymbol{x}_1$ and $\boldsymbol{x}_2$ are near some fixed location $\boldsymbol{x}_0 \in \mathbb{S}^2$, and $b_n(\cdot)$ is a sufficiently smooth function, for each $n\in\mathbb{N}_0$, we have  
$$  C(\boldsymbol{x}_1,\boldsymbol{x}_2) \approx  \sum_{n=0}^\infty  b_n(\boldsymbol{x}_0) {P}_n(\boldsymbol{x}_1^\top \boldsymbol{x}_2), \qquad \boldsymbol{x}_1,\boldsymbol{x}_2\in\mathbb{S}^2.$$
As a result, covariance functions with arbitrary locally isotropic behaviours are possible, making this class of models attractive.  In particular, the fractal index can vary from place to place on the spherical surface, in which case we say that the random field is multifractal.

\begin{ejemplo}
A natural extension of the Legendre-Mat\'ern model (\ref{guinness}) is given by the adaptive Schoenberg sequence
\begin{equation}
\label{spectrum}
 b_n(\boldsymbol{x})  = (\tau^2 + n^2)^{-\nu(\boldsymbol{x})-1/2}, \qquad n\in\mathbb{N}_0, \quad \boldsymbol{x}\in\mathbb{S}^2,
\end{equation}
where $\nu(\cdot)$ is a positive function that controls the local regularity of the sample paths.  In what follows, we call (\ref{spectrum}) the adaptive Legendre-Mat\'ern model. Of course, $\tau$ may also be a locally adaptive parameter, however, we essentially focus on the fractal properties of the model.
\end{ejemplo}

\subsection{Connections with a kernel-based method}

We now show that the proposed model, based on adaptive Schoenberg sequences, offers great flexibility. Actually, it is capable of emulating covariance functions coming from kernel-based methods, widely used in the spatial analysis literature \citep{fuentes2002spectral, nott2002estimation}. 

In a kernel-based approach, the random field $Z(\boldsymbol{x})$ on $\mathbb{S}^2$ is represented as a spatially weighted combination of isotropic random fields. This strategy permits to model dissimilar local dependency structures in different spatial zones.  Let  $D_1,\hdots,D_m$ be a collection of subregions that cover $\mathbb{S}^2$, and let $\lambda_j(\boldsymbol{x})$ be a positive kernel function centered at the centroid of $D_j$, for all $j=1,\hdots,m$ (see \citealp{schreiner1997locally}). 
Consider the random field
 \begin{equation*}
\label{fuentes}
Z(\boldsymbol{x}) = \sum_{j=1}^m \lambda_j(\boldsymbol{x}) Z_j(\boldsymbol{x}), \qquad \boldsymbol{x}\in\mathbb{S}^2,
\end{equation*} 
where $Z_1(\boldsymbol{x}),\hdots, Z_m(\boldsymbol{x})$ is a collection of independent isotropic random fields on $\mathbb{S}^2$.  Suppose that the covariance function of $Z_j(\boldsymbol{x})$ has isotropic part $K_j$. Thus, $Z(\boldsymbol{x})$ has the following covariance function
\begin{equation}
\label{cov_fuentes}
C(\boldsymbol{x}_1,\boldsymbol{x}_2) = \sum_{j=1}^m \lambda_j(\boldsymbol{x}_1) \lambda_j(\boldsymbol{x}_2) K_j\{d(\boldsymbol{x}_1,\boldsymbol{x}_2)\}, \qquad \boldsymbol{x}_1, \boldsymbol{x}_2\in\mathbb{S}^2.
\end{equation}

Since $K_j$ admits a representation of the form $K_j(d) = \sum_{n=0}^\infty p_{n,j} P_n (\cos d)$,  for each $j = 1,\hdots,m$, where $\{p_{n,j}: n\in\mathbb{N}_0\}$ is the associated $j$th Schoenberg sequence, the covariance function (\ref{cov_fuentes}) can be written as
\begin{equation}
\label{cov_fuentes2}
C(\boldsymbol{x}_1,\boldsymbol{x}_2) =\sum_{n=0}^\infty    \left\{ \sum_{j=1}^m \lambda_j(\boldsymbol{x}_1) \lambda_j(\boldsymbol{x}_2) p_{n,j} \right\} P_n(\boldsymbol{x}_1^\top \boldsymbol{x}_2),  \qquad \boldsymbol{x}_1, \boldsymbol{x}_2\in\mathbb{S}^2.
\end{equation}
We observe that (\ref{cov_fuentes2}) can be obtained as the sum of $m$ covariance functions of the form (\ref{cov_nonstationary}), where the $j$th adaptive Schoenberg sequence of type (\ref{cov_nonstationary}) can be taken as $b_{n,j}(\boldsymbol{x}) = \lambda_j^2(\boldsymbol{x}) p_{n,j}$.

\subsection{Simulation algorithm}
\label{section-sim}

This section presents a method for simulating anisotropic Gaussian random fields on $\mathbb{S}^2$. The representation (\ref{cov_nonstationary}) allows for an immediate simulation procedure based on the adaptive Schoenberg sequence. The following proposition is crucial to develop the simulation algorithm. 

\begin{proposition}
\label{prop1}
 Let $\kappa$ be a discrete random variable with ${\rm pr}(\kappa = n) = a_n$,  $n\in\mathbb{N}_0$, where $\{a_n: n\in\mathbb{N}_0\}$ is a probability mass sequence with a support containing that of the adaptive Schoenberg sequence $\{b_n(\boldsymbol{x}): n\in\mathbb{N}_0\}$, for all $\boldsymbol{x}\in\mathbb{S}^2$,  and let $\boldsymbol{w}$ be a random vector uniformly distributed on $\mathbb{S}^2$. Suppose that $\kappa$ and $\boldsymbol{w}$ are independent. Then, the random field
\begin{equation}
\label{simulacion2}
{Z}(\boldsymbol{x)} = \left\{ \frac{b_\kappa(\boldsymbol{x}) (2\kappa + 1)}{a_\kappa } \right\}^{1/2}    {P}_\kappa(\boldsymbol{w}^\top \boldsymbol{x}), \qquad \boldsymbol{x}\in\mathbb{S}^2,
\end{equation}
has mean $\mu(\boldsymbol{x}) = \{a_0 b_0(\boldsymbol{x})\}^{1/2}$, and its covariance function is given by (\ref{cov_nonstationary}).
\end{proposition}
This method appears as the spherical counterpart of the \emph{spectral turning bands} simulation algorithm used in Euclidean spaces (see, e.g.,  \citealp{mantoglou1982turning} and \citealp{emery2018continuous}).   For a neater exposition, the proof of Proposition \ref{prop1} is deferred to Appendix B.

  Even though the random field ${Z}(\boldsymbol{x})$ in (\ref{simulacion2}) has the predefined covariance function, it is clearly non-Gaussian distributed. By a central limit effect, an approximately Gaussian random field can be obtained by
\begin{equation*}
\widetilde{Z}(\boldsymbol{x}) =   \frac{1}{L^{1/2}}\sum_{\ell = 1}^L {Z}_\ell(\boldsymbol{x}), \qquad \boldsymbol{x}\in\mathbb{S}^2,
\end{equation*}
where ${Z}_1(\boldsymbol{x}),\hdots, {Z}_L(\boldsymbol{x})$ are $L$ independent copies from (\ref{simulacion2}), and $L$ must be a large integer.  

Basically, the method consists of an adequate weighted and rescaled combination of Legendre waves. The algorithm separates the choice of the adaptive Schoenberg sequence from the choice of the sequence $\{a_n: n\in\mathbb{N}_0\}$, which is equivalent to an important sampling technique. We observe that location and dispersion parameters can also be added in order to control the mean and variance of the sample paths. As discussed in \cite{emery2018continuous}, the process time of this algorithm is proportional to the number $L$ of copies and to the number of target locations, and it turns to be considerably fast.

\section{Numerical Examples}
\label{numerical}

\subsection{Simulated particles}
\label{examples}

\begin{figure}
\centering
 \includegraphics[scale=0.6]{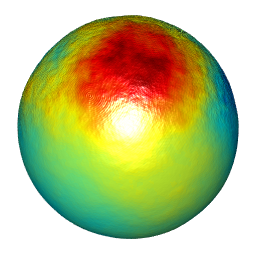}  \hspace{1cm}  \includegraphics[scale=0.6]{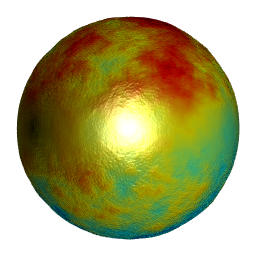}
\caption{Graphs showing the realizations of two star-shaped particles with locally varying fractal index according to scenarios (A) (left) and (B) (right).}
\label{ejemplo1}
\end{figure}

We now illustrate simulated particles from our proposal.   Let $\varphi_x \in [0,\pi]$ be the polar angle of $\boldsymbol{x}\in\mathbb{S}^2$ in the spherical coordinates system. We consider the adaptive Legendre-Mat\'ern model  (\ref{spectrum}), with $\tau^2=0.1$ and two different structures for $\nu(\boldsymbol{x})$: 
\begin{description}
\item[(A)] Increasing Hausdorff dimension, from South to North, according to $$\nu(\boldsymbol{x}) = 0.6 + 0.3 \times [1+\exp\{-10 (\varphi_x-\pi/2)\}]^{-1}.$$
\item[(B)] Hausdorff dimension with high values at the poles and lower values at the Equator, according to $$\nu(\boldsymbol{x}) = 0.6 + 0.3 \times \exp\left\{-20 (\varphi_x-\pi/2)^2  \right\}.$$ 
\end{description}
 For each case, the Hausdorff dimension varies in terms of the polar angle, approximately in the range $(2.1, 2.4)$.  We set $b_0(\boldsymbol{x}) \equiv 0$ and incorporate a constant mean. We have also rescaled the covariance function to obtain a constant variance. The resulting radial function has mean $E\{Z(\boldsymbol{x})\} = 100$ and variance ${\rm var}\{Z(\boldsymbol{x})\} = 10$.  The random variable $\kappa$ in Proposition \ref{prop1} is chosen to follow a negative binomial distribution with parameters $r = 1$ (number of failures) and $p = 0.1$ (success probability). For this distribution, the probability mass function decays slowly to zero, i.e. it favors the occurrence of large values of $\kappa$,  which is a desirable property in the sampling process. In our experiment, we have considered a spatial grid, in terms of the polar and the azimuthal angles, with $2000 \times 1000$ target points. Figure \ref{ejemplo1} shows the simulated particles, with $L=300$, for the cases (A) and (B), where it is seen that the simulations match the theoretical features. 

\subsection{Earth topography: continents versus seafloor roughness}

The aim of this section is to provide a simulated version of the Earth topography. We consider an extension of the monofractal simulation given by \cite{hansen2015gaussian}.  In our analysis, we take into account the different degrees of roughness present in the continents and the seafloor \citep{gagnon2006multifractal}. 

The mean function $\mu(\boldsymbol{x})$ is estimated from online data sources (see, e.g., data outputs from the Joint Institute for the Study of the Atmosphere and Ocean, Seattle, United States), by using spherical harmonics regression, which is the natural basis for the spherical geometry \citep{marinucci2011random}, i.e. we consider the regressors $P_n^m(\sin \varphi_x) \cos(m \theta_x)$ and $P_n^m(\sin \varphi_x) \sin(m \theta_x)$, where  $\varphi_x$ and $\theta_x$ are the polar and azimuthal angles in the spherical coordinates system, and $P_n^m$ denotes the associated Legendre polynomial, for $n\in\mathbb{N}_0$ and $m = 0, \hdots, n$. We set $n=60$ in the estimation setting.

\begin{figure}
\centering
 \includegraphics[scale=0.45]{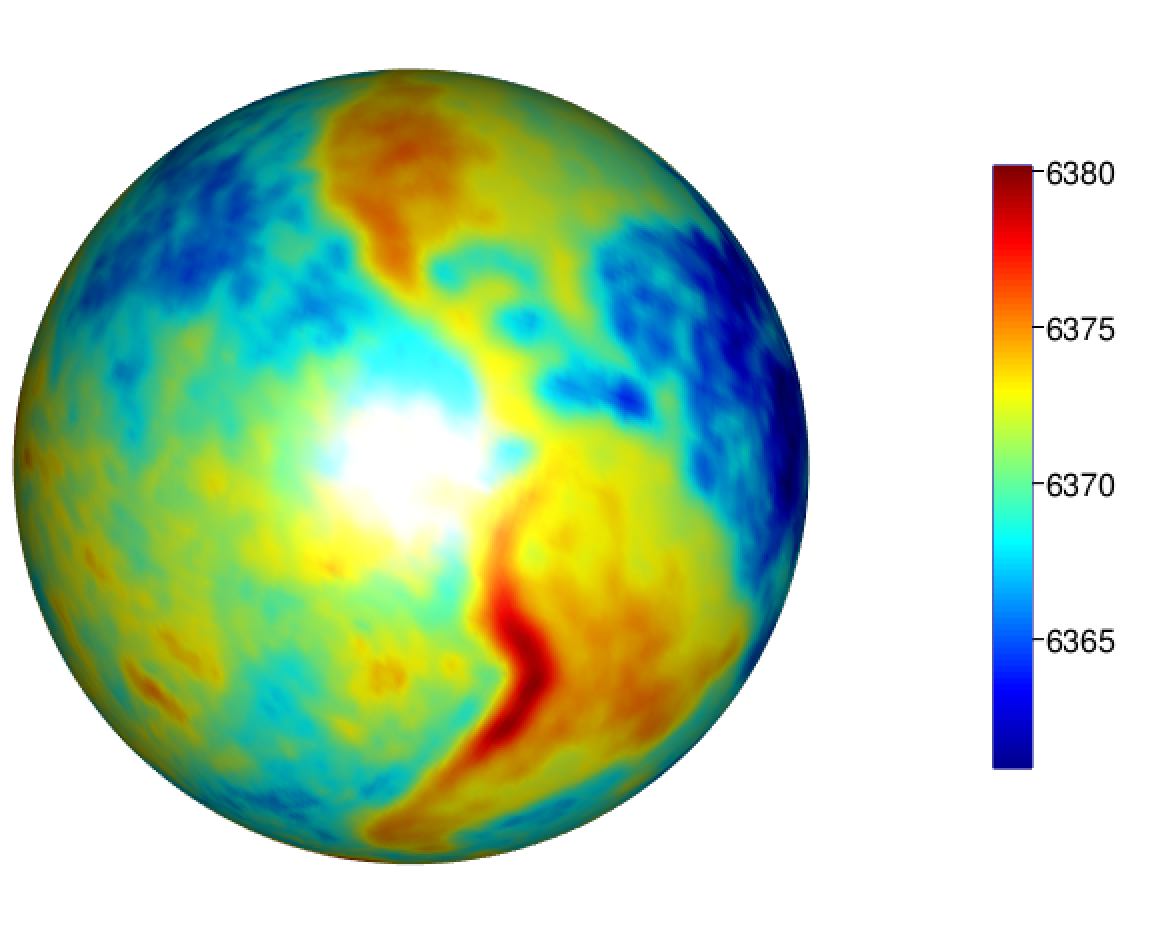} 
\caption{Simulated Earth surface, where the elevation is measured in kilometers. The Hausdorff dimensions considered are $\eta_{0,S} = 2.68$ for the seafloor and $\eta_{0,C} = 2.46$ for  continents.}
\label{planeta}
\end{figure}

While the mean function captures the large scale variation of the Earth topography, we must incorporate the small scale variation.  According to \cite{gagnon2006multifractal}, we distinguish between the Hausdorff dimension of the seafloor ($\eta_{0,S} = 2.68$) and the continents ($\eta_{0,C} = 2.46$). We then consider a dichotomic Schoenberg sequence from the adaptive Legendre-Mat\'ern model. We rescale the covariance, and the resulting realization is chosen to have a constant standard deviation of approximately $2640$ meters, obtained empirically from the online resources mentioned above.  Figure \ref{planeta}  depicts the simulated topography of planet Earth, on a grid of longitudes and latitudes with resolution of $1^\circ\times 1^\circ$, where we have used the same setting as in Section \ref{examples}.

\begin{figure}
\centering
 \includegraphics[scale=0.5]{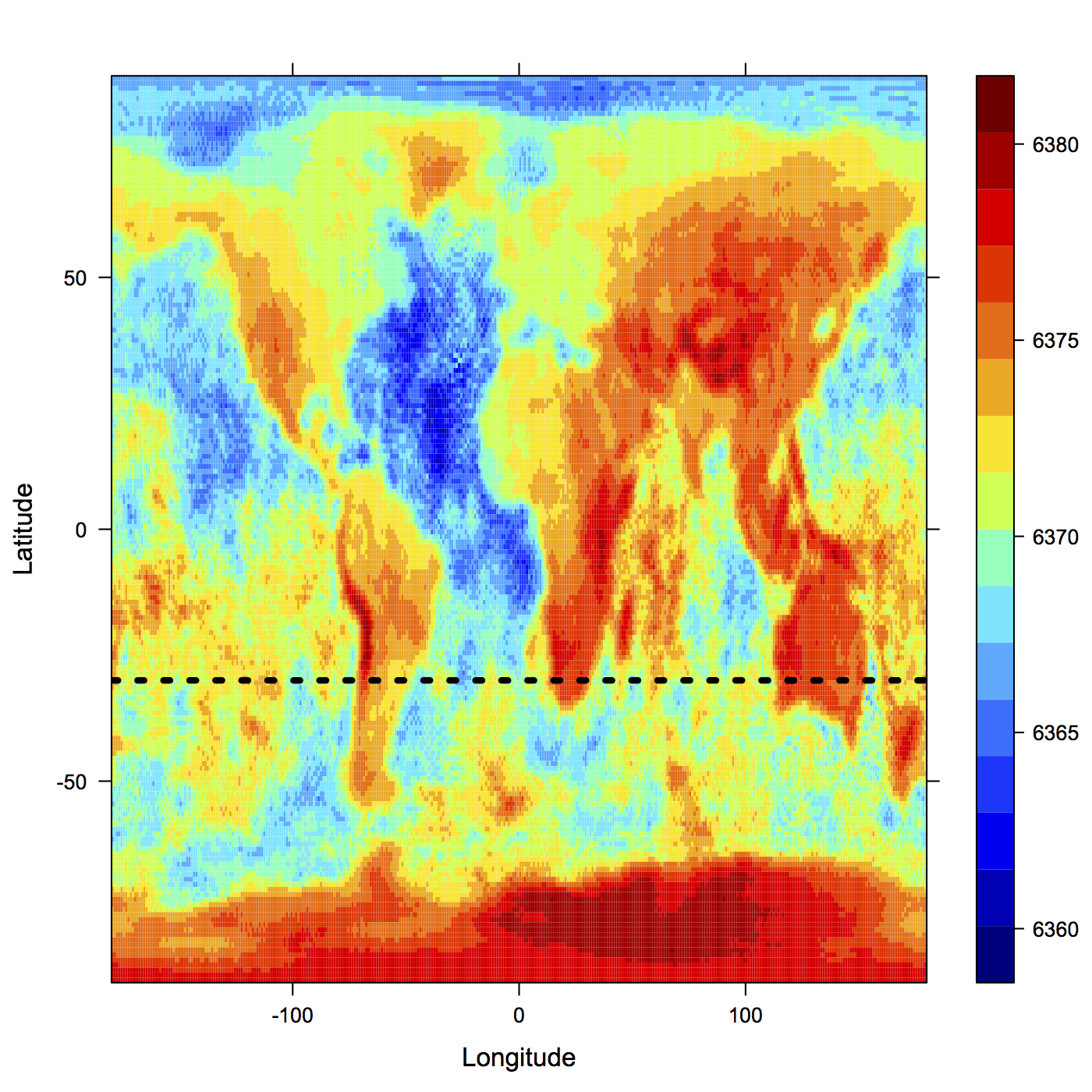}  \includegraphics[scale=0.5]{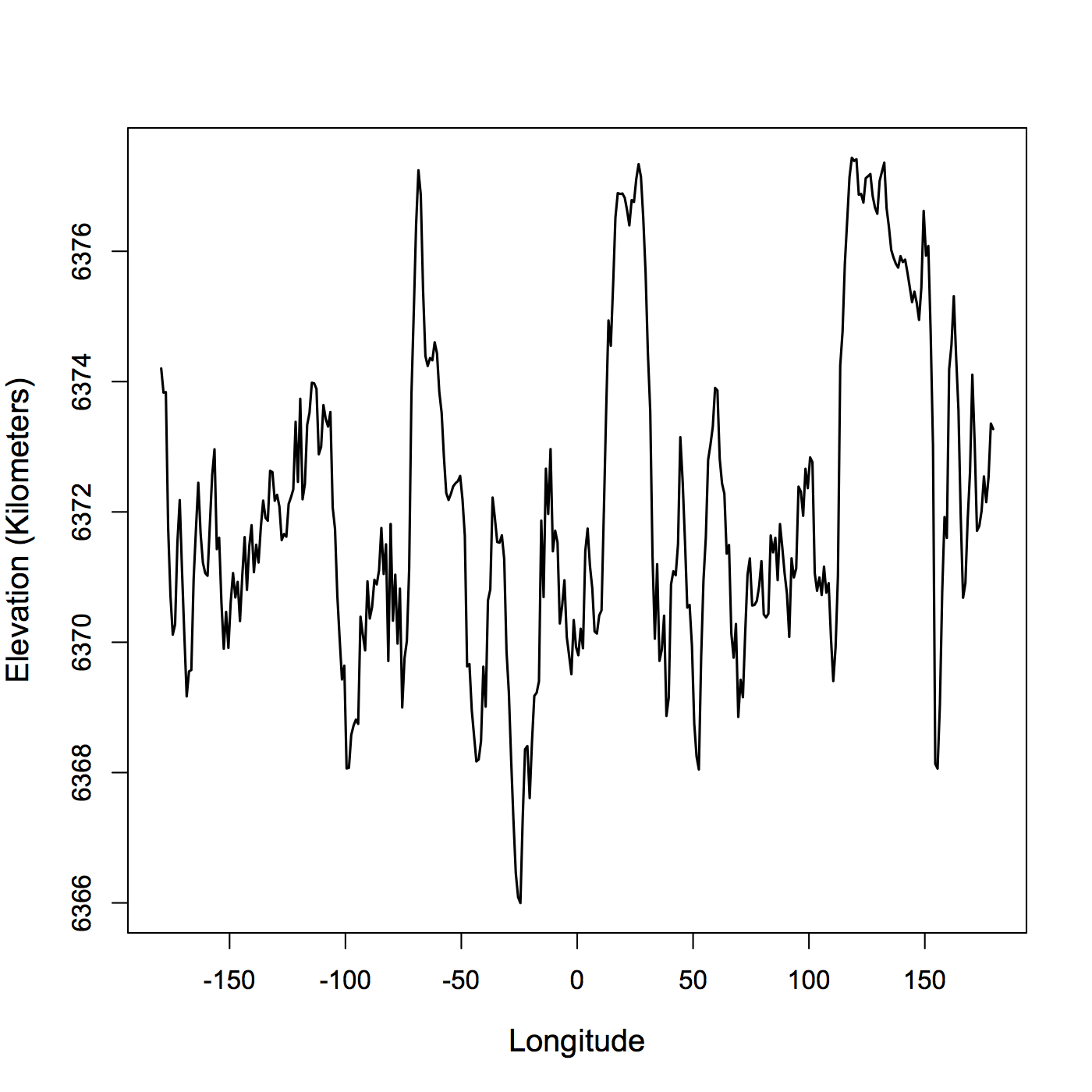}
\caption{Graphs showing the simulated Earth topography on a grid of longitudes and latitudes (top) and the radial function along the latitude $-30^\circ$ (bottom). The elevation is measured in kilometers.}
\label{elevation}
\end{figure}

Figure \ref{elevation} displays the Earth topography as a function of the longitude and latitude. It gives us a complete picture of the Earth map. The radial function (in kilometers) along the latitude $-30^\circ$ is also supplied. This example illustrates that continents are smoother than the seafloor, which is justified by tectonic forces and erosion-like processes \citep{gagnon2006multifractal}.

\section{Discussion}

\label{discusion}

We have introduced a flexible framework for modelling and simulating three-dimensional multifractal star-shaped particles. The radial function of the object has been modelled by mean of an anisotropic Gaussian random field on the sphere, which is generated from a locally varying Schoenberg sequence. A simple adaptive version of the Legendre-Mat\'ern covariance function has been proposed and employed to generate particles with place to place variable Hausdorff dimension. Our findings have been exemplified through numerical experiments, including an illustration of the Earth topography. Our approach may be used as a building block for more sophisticated models emulating the surface of celestial bodies. 

 A natural generalization of this work is to consider temporally dependent Schoenberg sequences, allowing for particles with a dynamically updated Hausdorff dimension. The findings of \cite{berg2017schoenberg} might be useful here. Another interesting problem is the search for efficient methods for estimating locally variable Hausdorff dimensions. One challenging possibility is to adapt the tools reviewed in \cite{gneiting2012estimators}.  Advances in this direction also include the works of \cite{anderes2011local} and \cite{ziegel2013stereological}. The extension to non-Gaussian particles can be tackled by using transformations of Gaussian random fields as in \cite{xu2017tukey}. This  topic is interesting because the probability distribution of the random field can influence some geometrical aspects of the particle. For instance, \cite{hansen2015gaussian} show that Gamma-L\'evy particles exhibit more pronounced spikes. From a computational viewpoint, one may be interested in proposing an  improved version of the simulation algorithm by using parallel computing.

Finally, \cite{gneiting2013strictly} states, in his open Problem 15, that new methodologies involving anisotropic dependencies are also desirable in environmental and climatological phenomena (see also \citealp{hitczenko2012some} and \citealp{castruccio2013global}). So we believe that our findings in Section \ref{proposal} can be useful to develop new applications in various fields related to spatial analysis.



\appendix
\section*{Appendices}
\addcontentsline{toc}{section}{Appendices}
\renewcommand{\thesubsection}{\Alph{subsection}}
\renewcommand{\theequation}{\Alph{subsection}.\arabic{equation}}
\renewcommand\thetable{\thesubsection.\arabic{table}}

\counterwithin{thm}{subsection}
\counterwithin{mydef}{subsection}

\subsection{Positive definiteness of (\ref{cov_nonstationary})}
\label{apendice2}

The set of spherical harmonic functions,  $\{{Y}_{n,m}: n\in\mathbb{N}_0, m = -n, \hdots, n\}$,  form an orthogonal basis of the Hilbert space of square integrable functions on $\mathbb{S}^2$. Explicit expressions for the real and the imaginary parts of $Y_{m,m}$ has been used in Section \ref{numerical}. The \emph{addition theorem} for spherical harmonic functions \citep{marinucci2011random} establishes that 
\begin{equation*}
  {P}_n(\boldsymbol{x}_1^\top \boldsymbol{x}_2) = \sum_{m=-n}^n {Y}_{n,m}(\boldsymbol{x}_1)   \overline{ {Y}_{n,m}(\boldsymbol{x}_2) }, \qquad \boldsymbol{x}_1, \boldsymbol{x}_2\in\mathbb{S}^2,
\end{equation*}
where $\overline{c}$ denotes the complex conjugate of $c\in\mathbb{C}$. The semi positive definiteness of (\ref{cov_nonstationary}) is a direct consequence of the addition theorem. In fact, a straightforward calculation shows that
\begin{eqnarray*}
\sum_{i=1}^k \sum_{j=1}^k  a_i a_j C(\boldsymbol{x}_i,\boldsymbol{x}_j)    &  =   & \sum_{i=1}^k \sum_{j=1}^k  a_i a_j  \sum_{n=0}^\infty   \left\{ b_n(\boldsymbol{x}_i)  b_n(\boldsymbol{x}_j) \right\}^{1/2}   {P}_n(\boldsymbol{x}_i^\top \boldsymbol{x}_j)\\
&  =   & \sum_{i=1}^k \sum_{j=1}^k  a_i a_j  \sum_{n=0}^\infty   \left\{ b_n(\boldsymbol{x}_i)  b_n(\boldsymbol{x}_j) \right\}^{1/2}   \sum_{m=-n}^n {Y}_{n,m}(\boldsymbol{x}_i)  \overline{ {Y}_{n,m}(\boldsymbol{x}_j) }\\
&  =   &  \sum_{n=0}^\infty   \sum_{m=-n}^n   \left| \sum_{i=1}^k   a_i   \{ b_n(\boldsymbol{x}_i)  \}^{1/2}   {Y}_{n,m}(\boldsymbol{x}_i) \right|^2,
\end{eqnarray*}
where $|c| = (\overline{c} c)^{1/2}$ denotes the magnitude of $c$. The last expression is clearly nonnegative.

\subsection{Proof of Proposition \ref{prop1}}
\label{apendice1}

Before we state the proof of Proposition \ref{prop1}, we must introduce some properties of Legendre polynomials.  \citet[Equation (2.4)]{ziegel2014convolution} establishes that, for any $n,k\in\mathbb{N}_0$,
\begin{equation}
\label{propiedad1}
\int_{\mathbb{S}^2}   {P}_n(\boldsymbol{w}^\top \boldsymbol{x}_1) {P}_k(\boldsymbol{w}^\top \boldsymbol{x}_2)  U({\rm d}\boldsymbol{w})   =  \frac{\delta_{n,k}}{2n+1} {P}_n(\boldsymbol{x}_1^\top \boldsymbol{x}_2),  \qquad \boldsymbol{x}_1,\boldsymbol{x}_2\in\mathbb{S}^2,
\end{equation}
where $U$ is the uniform probability measure on $\mathbb{S}^2$ and $\delta_{n,k}$ denotes the Kronecker delta.  In particular, for any $n\in\mathbb{N}$, we have 
\begin{equation}
\label{propiedad2}
\int_{\mathbb{S}^2}   {P}_n(\boldsymbol{w}^\top \boldsymbol{x})  U({\rm d}\boldsymbol{w})   = 0,  \qquad \boldsymbol{x}\in\mathbb{S}^2.
\end{equation}

We now proceed with the proof of Proposition \ref{prop1}. Let ${Z}(\boldsymbol{x})$ be the random field defined in (\ref{simulacion2}). The mean function is 
\begin{equation*}
E\{{Z}(\boldsymbol{x})\}    =     \sum_{n=0}^\infty   a_n  \left\{ \frac{b_n(\boldsymbol{x}) (2n + 1)}{a_n} \right\}^{1/2}  \int_{\mathbb{S}^2}  {P}_n(\boldsymbol{w}^\top \boldsymbol{x})   U(\text{d}{\boldsymbol{\omega}}).
 \end{equation*}
The property (\ref{propiedad2}) implies that $E\{{Z}(\boldsymbol{x})\} = \{a_0 b_0(\boldsymbol{x})\}^{1/2}$. On the other hand, the covariance function of ${Z}(\boldsymbol{x})$ is 
\begin{equation*}
{\rm cov}\{{Z}(\boldsymbol{x}_1), {Z}(\boldsymbol{x}_2)\}     =     \sum_{n=0}^\infty \{b_n(\boldsymbol{x}_1) b_n(\boldsymbol{x}_2)\}^{1/2} (2n+1) \int_{\mathbb{S}^2}  {P}_n(\boldsymbol{w}^\top \boldsymbol{x}_1) {P}_n(\boldsymbol{w}^\top \boldsymbol{x}_2) U(\text{d}{\boldsymbol{\omega}}), \qquad \boldsymbol{x}_1,\boldsymbol{x}_2\in\mathbb{S}^2.
 \end{equation*}
Using (\ref{propiedad1}), we obtain the covariance function (\ref{cov_nonstationary}). The proof is completed.

\bibliographystyle{apalike}
\bibliography{mybib}

\end{document}